\documentclass[12pt,a4paper]{amsart}

\usepackage{amsmath,eucal}
\usepackage{amssymb}

\newtheorem{thm}{Theorem}[section]

\newtheorem{lem}[thm]{Lemma}
\newtheorem{prp}[thm]{Proposition}
\newtheorem{cor}[thm]{Corollary}

\theoremstyle{definition}

\newtheorem{rmk}[thm]{Remark}

\newcommand{\la}{\langle}
\newcommand{\ra}{\rangle}

\title{Conformality of a differential with respect to
Cheeger-Gromoll type metrics}
\author{W. Koz\l owski \and K. Niedzia\l omski}

\begin{document}

\begin{abstract} We investigate conformality of the differential of a mapping between Riemannian manifolds if the tangent bundles are equipped with
a generalized metric of Cheeger-Gromoll type.
\end{abstract} 


\maketitle

\section{Introduction and Preliminaries}
Generalized metrics of Cheeger-Gromoll type or $(p,q)$-metrics $h_{p,q}$, being a generalization of Sasaki metric $h_S$ \cite{Kw}
and Cheeger-Gromoll metric  $h_{CG}$ \cite{GK}, have been recently introduced by M. Benyounes, E. Loubeau and C. M. Wood 
in \cite{BLW}  in the context of harmonic sections. In \cite{BLW1}, 
the same authors showed that the geometry of the tangent bundle equipped with this kind of metric is of independent interest. It is worth noticing that M. I. Munteanu in \cite{M} investigated independently 
the geometry of tangent bundle equipped with a certain deformation of Cheeger-Gromoll metric other that in \cite{BLW}. Yet in  \cite{KW}, Sz. Walczak 
and the first named author considered
$(p,q)$-metrics in the context of Riemannian submersions and Gromov-Hausdorff topology.

In this paper we introduced $(p,q,\alpha)$-metrics which are more general that $(p,q)$-metrics. (In contrast to \cite{BLW} we do not
assume that $p,q$ and $\alpha$ are constant).
We investigate relations between conformality of a map $\varphi:(M,g)\to (M',g')$ between Riemannian manifolds and its differential $\Phi=\varphi_\ast: (TM, h)\to 
(TM',h')$
 between their tangent bundles equipped with $(p,q,\alpha)$-metric $h$ and $(r,s,\beta)$-metric $h'$, respectively.

Interesting enough, there is essential difference between the cases $\dim M=2$ and $\dim M\geq 3$. 

We prove that in the second case
(Theorem \ref{MR_T1})  $\Phi$ is conformal if and only if $\varphi$ is a homothety and totally geodesic immersion
 and some special relations between triples $(p,q,\alpha)$ and $(r,s,\beta)$ hold. In this case $\Phi$ is also a homothety with the same dilatation as
 $\varphi$.

However, in the first case it may happen that $\Phi $ is conformal, 
although $\varphi$ is not a totally geodesic immersion (Theorem \ref{MR_T2}). Then $\Phi$ is  no longer a homothety. An example of such a map is 
given.

\subsection{Cheeger-Gromoll type metrics}\label{S_CGM}
Consider a Riemannian manifold $(M,g)$, and let $\pi :TM \to M$ be its tangent bundle.  
The Levi-Civita connection $\nabla$ of $g$, gives a natural splitting $T(TM) = \mathcal H \oplus \mathcal V$
of the second tangent bundle $\pi_\ast :T(TM)\to TM$, where the {\it vertical} distribution $\mathcal V$ is the kernel of $\pi_\ast$, and the {\it 
horizontal} distribution is the kernel of, so called, connection map $K$. 
If $X,Z\in T_x M$ then by $X^v_Z$ we denote the vertical lift of $X$ to the point $Z$, i.e., $X^v_Z$ is a tangent vector to 
the curve $t \mapsto Z + tX$ at $t=0$. 
Every $A\in T_Z(TM)$ can be uniquely written as $A=\mathcal H A+ \mathcal V A$, where  
$\mathcal H A\in \mathcal H_Z$ 
and $\mathcal V A\in \mathcal V_Z$ denote its horizontal and vertical part respectively.
The vertical part of $A$ is given  by $(KA)^v_Z$.

Recall that $K$ is a smooth $\mathbb R$-linear bundle morphism 
determined by the conditions:
\begin{itemize}
\item[(K1)] For every 
$Z\in TM$, 
$K: T_Z(TM)\to T_{\pi(Z)} M$ 
is the canonical isomorphism, i.e., $K(X^v_Z)=X$.
\item[(K2)] For every vector field $X$ on $M$ and every $v\in T_xM$, $K(X_\ast v) = \nabla_v X$.
\end{itemize}   
Notice that (K1) and (K2) imply the following properties
 
\begin{itemize}
\item[(K3)] For every Riemannian manifold $(M', g')$ and every $X,Z\in T_x M$ and every map $\varphi: M\to M'$, $\varphi_{\ast \ast} X^v_Z  = 
(\varphi_\ast X)^v_{\varphi_\ast (Z)}$.
\item[(K4)] For every curve $\gamma$ in $M$ and every vector field $\xi$ along $\gamma$, $K(\dot \xi)=\nabla_{\dot \gamma} \xi$. 
\end{itemize}

Let $p,q,\alpha$ be smooth functions on $M$. Assume $q$ is non-negative and $\alpha$ is positive. Define $(p,q,\alpha)$-metric $h=h_{p,q,\alpha}$ on $TM$ 
as follows:
For every $A,B\in T_Z(TM)$, $Z\in T_x M$,
\[h(A,B) = g(\pi_\ast A,\pi_\ast B) +
\omega_{\alpha}(Z)^{p}
\big( g(KA,KB) + qg(KA,Z)g(KB,Z)\big),\]
where  $\omega_\alpha(Z)=(1+\alpha g(Z,Z))^{-1}$. Here all functions $p,q,\alpha$ are evaluated at $x$. 
For any $p,q,\alpha$, the Riemannian metric $h_{p,q,\alpha}$ is a special case of a metric considered in \cite{M}.
Notice that if $p,q,\alpha$ are constants and $\alpha =1$ then $h_{p,q,\alpha}$ becomes a metric from \cite{BLW}. In particular,  
$h_{0,0,1}$ (resp. $h_{1,1,1}$) is Sasaki metric $h_S$ \cite{Kw} (resp.
 Cheeger-Gromoll metric $h_{CG}$ \cite{GK}).\\


\subsection{Conformal mappings and metrics}\label{SS_CM} 
Recall that a map $\varphi: (M,g)\to (M',g')$ between Riemannian manifolds is {\it conformal} if $\varphi^\ast g' = \lambda g$ 
for some positive function  $\lambda$ on $M$. The function $\lambda$ is called  a {\it dilatation}. 
A conformal mapping with constant dilatation is called a {\it homothety}. If $\dim M < \dim M'$ a conformal mapping $\varphi$ is often called
{\it weakly conformal}. 

Let $\lambda$ be a strictly positive $C^\infty$-function on $M$, and $g^\lambda= \lambda g$. The Levi-Civita connections 
$\nabla$ and $\nabla^\lambda$ of $g$ and $g^\lambda$ are related as follows: $\nabla^\lambda = \nabla + S_M^{g,\lambda}$ where $S=S_M^{g,\lambda}$
is a symmetric $(1,2)$-tensor field given by (compare \cite{Ch}, page 64):
\[S(X,Y) = \frac{1}{2\lambda}\big( (X\lambda)Y + (Y\lambda)X - g(X,Y) {\rm grad\,} \lambda \big),
\quad X,Y \in \Gamma(M, TM).\]

Suppose $\varphi:M\to M'$ is an immersion, e.g., conformal mapping. Then for every $x\in M$ we may choose an open neighbourhood
$U_x$ of $x$ such that $L_{x'}=\varphi(U_x)$ ($x'=\varphi(x)$) is a regular submanifold of $M'$.  Let $j:L_{x'}\to M'$ be the inclusion map, and let
$\bar g =j^\ast g'$ be the induced metric tensor on $L_{x'}$. 
Moreover, let $\Pi$ denote
the second
fundamental form of $L_{x'}$.
We say that the immersion $\varphi:(M,g)\to (M',g')$ is {\em totally geodesic} if
for every $x\in M$, $L_{x'}$ is a totally geodesic submanifold of $(M',g')$. One can prove the following

\begin{lem}\label{CM_nable} Suppose $\varphi: (M,g)\to (M',g')$ is a conformal mapping with a dilatation $\lambda$. Choose $x\in M$ and put $x'=\varphi(x)$.
Let $\gamma:(-\epsilon,\epsilon)\to L_{x'}$ be a curve in $M$ and $\xi$
be a vector field along $\gamma$. Put $\gamma'=\varphi\circ \gamma$ and $\xi'=\varphi_\ast \xi$. Then
\[ \bar \nabla_{\dot \gamma'}\xi' = \varphi_\ast\nabla_{\dot \gamma}\xi + \varphi_\ast S(\dot \gamma, \xi),\]
where $S=S_M^{g,\lambda}$, and $\nabla$ and $\bar \nabla$ are the Levi-Civita connections of $g$ and $\bar g$ respectively.
\end{lem}

Adopt the notations from Lemma \ref{CM_nable}. Put $Z= \xi(0)$, $Z'=\xi'(0)$, $v= \dot\gamma (0)$ and $v'=\dot\gamma'(0)$. Suppose that a vector $A\in 
T_Z(TM)$ is tangent to the curve $\xi$
(it is convenient to think of vector fields along curves as of curves in the tangent bundle), i.e., $A=\dot \xi(0)$. Next, let $K$ and 
$K'$ denote connection maps induced from $\nabla$ and $\nabla'$ respectively. Moreover put $\Phi=\varphi_\ast:TM\to TM'$.
As a direct consequence of Lemma \ref{CM_nable}, the equation $\nabla'=\nabla+\Pi$ and properties of connection map we get

\begin{lem}\label{CM_KK} The vectors $K(A)$ and $K'(\Phi_\ast A)$ are related as follows:
\[ K'(\Phi_\ast A) = \varphi_\ast K(A) + \varphi_\ast S(v,Z)+\Pi(v',Z').\]
In particular, if $A$ is horizontal then $ K'(\Phi_\ast A) = \varphi_\ast S(v,Z)+\Pi(v',Z').$
\end{lem}

\begin{rmk} If $\dim M=\dim M'$ then the term $\Pi$ is omitted.
\end{rmk}

Let $U=U_x$ and $L'=L_{x'}$ and let $\pi':TL'\to L'$ be a natural projection.
Since $\varphi: U\to L'$ is a conformal diffeomorphism, so is its inverse. In particular $\Phi: TU\to TL'$ is a diffeomorphism. Therefore we have
\begin{cor} \label{CM_cor_-1}
Take $A'\in T_{Z'}(TL')$ such that $v'=\pi'_\ast A$. Then
\[ K(\Phi_\ast^{-1} A')= \varphi_\ast^{-1} K'(A') + \varphi_\ast^{-1} S'(v',Z')-\varphi_\ast^{-1}\Pi(v',Z'),\]
or equivalently
\[ K(\Phi_\ast^{-1} A')= \varphi_\ast^{-1} \bar K(A') + \varphi_\ast^{-1} S'(v',Z'),\]
where $S'=S_{L'}^{\bar g, \mu}$ with $\mu= (1/\lambda)\circ \varphi^{-1}$, and $\bar K$ is the connection map induced from $\bar \nabla$.
\end{cor}

\begin{cor}\label{CM_co} Suppose $\varphi:(M,g)\to (M',g')$ is a
conformal mapping and  $M$ is connected. $\Phi_\ast$ maps horizontal vectors onto horizontal vectors if and only if $\varphi$ is a totally geodesic homothety.
\end{cor}

\begin{proof} ($\Rightarrow$) If $\Phi_\ast$ maps horizontal vectors onto horizontal vectors then by Lemma \ref{CM_KK},
$\varphi_\ast S + \Pi$ vanishes identically. Since $\varphi_\ast S$ and $\Pi$ are always orthogonal and a
conformal mapping is an immersion it follows that $S$ and $\Pi$ vanish identically. 
Applying the definition of $S$ with $X=Y={\rm grad\,} \lambda$, we get
that
${\rm grad\,} \lambda$ is the zero vector field. Consequently, $\lambda$ is constant and therefore $\varphi$ is a homothety. Since $\Pi$ vanishes,
$\varphi$ is totally geodesic.

($\Leftarrow$) Obvious.
\end{proof}

\subsection{ Algebraic lemmas}

Suppose two finite dimensional spaces $V$ and $W$ equipped with inner products 
$\la ,\ra_V $ and $\la ,\ra_W $ are given. 
Let $B:V\times V\to W$ be  a symmetric, bilinear form on $V$. Moreover, let $C\geq 0$. Consider a condition
\begin{equation}\label{SS_AL_E1}
\la B(X,Z),B(Y,Z)\ra_W = C\la X,Y\ra_V \la Z,Z\ra_V.
\end{equation}
for every $X,Y,Z\in V$. Then, if $X,Y$ are orthogonal
\begin{equation}\label{SS_AL_E01}
\la B(X,X),B(Y,Y)\ra_W =
 -C\la X,X\ra_V \la Y,Y\ra_V.
\end{equation}
\begin{lem}\label{Lem_AL}
Assume that $B$ satisfies the condition \eqref{SS_AL_E1}.
If $\dim V \geq 3$ then $C=0$. In particular, $B$ vanishes.
\end{lem}
\begin{proof}Suppose that $C\ne 0$.
Take an orthonormal  pair $X,Y$.  
 Let $\xi = B(X,X)$ and $\zeta = B(Y,Y)$.
 By \eqref{SS_AL_E1} we have
\begin{equation}\label{SS_AL_E2}
\la \xi,\xi\ra= \la \zeta,\zeta\ra = C>0. 
\end{equation}
In particular, $\xi\ne 0$ and $\zeta\ne 0$.
Applying \eqref{SS_AL_E01} we see that
\[\la \xi, \zeta \ra = - C. \]
Using above and \eqref{SS_AL_E2} one can obtain that
$\xi=-\zeta$. 
Next, since $\dim V\geq 3$ we may find $Z\in V$ such that $X,Y,Z$ is an orthonormal triple.
Let $\eta =  B(Z,Z)$. Then, by above $\xi=-\zeta=\eta=-\xi$, which contradicts the fact that $\xi\ne 0$.
\end{proof}
 Notice that the assumption $\dim V\geq 3$ is essential. Namely we have
\begin{lem}\label{SS_AL_PRP}
$($a$)$ A symmetric bilinear form 
$B:\mathbb{R}^2\times \mathbb{R}^2\to \mathbb{R}^2$ satisfies
\eqref{SS_AL_E1} if and only if there exists an angle $\theta$ such that
\begin{equation}\label{AL_E}
B(X,Y)=\pm\sqrt{C}e^{i\theta}XY\quad {\rm or}\quad B(X,Y)=\pm\sqrt{C}e^{i\theta}\bar{X}\bar{Y},
\end{equation}
where we identify $\mathbb{R}^2$ with $\mathbb{C}$. \\
$($b$)$ If $\dim V=2$ and a non-zero symmetric bilinear form $B:V\times V\to W$ satisfies the condition \eqref{SS_AL_E1} 
then there exists
a $2$-dimensional subspace $U$ of $W$ such that the image of $B$ is equal to $U$ and with respect to orthonormal bases of $V$ and $U$, $B$ is of the form 
\eqref{AL_E}.
\end{lem}

\begin{proof} (a) Elementary exercise. (b) Take an orthonormal basis $X,Y$ of $V$. Since $B\ne 0$, we have $C\ne 0$. Consequently, 
$\xi = B(X,X)$, $\zeta = B(Y,Y)$ and $\eta = B(X,Y)$ are nonzero vectors of $W$ of length  $\sqrt{C}$. 
By \eqref{SS_AL_E1} we see that $\la \xi,\eta\ra = \la  \zeta, \eta\ra = 0$.   Moreover,  $\la \xi, \zeta\ra = -C$, by  \eqref{SS_AL_E01}.
It follows that $\xi = -\zeta$. Consequently, the image $U$ of $B$ is a two-dimensional subspace spanned by $\xi,\eta$. 

Now taking orthonormal bases of $V$ and $U$, e.g., $X,Y$ and $\xi/\sqrt{C}, \eta/\sqrt{C}$, we reduce (b) to (a).
 \end{proof}

\section{Conformality of a differential}
In this section all manifolds are connected.
Let $(M,g)$ and $(M',g')$ be Riemannian manifolds of dimensions $m$ and $m'$, respectively.
We assume that $m,m'\geq 2$. Denote by $\nabla$ and $\nabla'$ the Levi-Civita connections of $g$ and $g'$, respectively.
Equip their tangent bundles $\pi:TM\to M$ and $\pi':TM'\to M'$
with $(p,q,\alpha)$-metric $h$ and $(r,s,\beta)$-metric $h'$, respectively.

Let $\varphi: M\to M'$ and $\Phi = \varphi_\ast: TM\to TM'$. Put $g=\la,\ra$ and $g'=\la,\ra'$.
Denote by $|\cdot|$ and $|\cdot|'$ the norms induced by $g$ and $g'$, respectively. Moreover, denote by $\|\cdot\|$ and $\|\cdot\|'$ the norms 
induced by $h$ and $h'$, respectively.

In the paper we use the following notation:
If $\varphi$ (resp. $\Phi$) is conformal mapping then its dilatation will be always denoted 
by $\lambda$ (resp. $\Lambda$).

\subsection{Technical lemmas}

\begin{lem}\label{C1} Suppose that $\varphi$ and $\Phi$ are conformal mappings. 
 Then  for any $Z\in T_x M$ and $x'=\varphi(x)$,
\begin{eqnarray}
\Lambda(Z) &=& \lambda(x) \frac{(1 + \alpha(x) |Z|^2)^{p(x)}}{(1 + \lambda(x) \beta(x') |Z|^2)^{r(x')}},\label{C11}\\ 
   q (x) &=& \lambda (x) s (x').\label{C12}
\end{eqnarray}
\end{lem}
\begin{proof} 
Let $X,Z\in T_xM$. Applying conformality of $\Phi$ and (K3) we have: 
$    \| (\varphi_\ast X)^v_{\varphi_\ast Z} \|'^2=\| \Phi_\ast X_Z^v \|'^2= \Lambda(Z) \|X^v_Z\|^2$. 
Using now the definitions of $h$ and $h'$ and conformality of $\varphi$, one can easily get
\begin{equation}\label{C13}
 \lambda(x) \frac{|X|^2 + \lambda(x) s(x') \la X,Z\ra^2}{(1 + \lambda(x) \beta(x') |Z|^2)^{r(x')}} = \Lambda(Z) \frac{|X|^2 + q(x) \la X,Z\ra^2}{(1+\alpha(x) 
|Z|^2)^{p(x)}}.
\end{equation} 

Taking nonzero vector $X$ orthogonal to $Z$, \eqref{C13} becomes \eqref{C11}. Next, let $Z\ne 0$. Putting $X=Z$ in \eqref{C13} and 
comparing the result with \eqref{C11} we get \eqref{C12}.
\end{proof}

\begin{lem}\label{C3}
If $\Phi$ is conformal then so is $\varphi$. Moreover $\lambda(x) = \Lambda(0_x)$, $x\in M$.
\end{lem}
\begin{proof}
Let $x\in M$ and  $Z=0_x\in T_xM$. By (K3) and conformality of $\Phi$, for every $X,Y\in T_xM$
\begin{eqnarray*}
\la\varphi_\ast X,\varphi_\ast Y\ra' &=& h'((\varphi_\ast X)^v_{\varphi_\ast Z},(\varphi_\ast Y)^v_{\varphi_\ast Z})\\
&=& \Lambda(Z) h(X^v_Z, Y^v_Z)=\Lambda(Z)\la X,Y\ra.
\end{eqnarray*}
\end{proof}

\begin{lem}\label{C2} 
Adopt the notation from Section \ref{SS_CM}.
Suppose that $\varphi$ and $\Phi$ are conformal mappings. 
Then 
\begin{itemize}
\item[(a)] $\varphi$ is a homothety.
\item[(b)]
For every $x\in M$ one of the following conditions holds:
\begin{eqnarray}
&&p(x)=r(x')=0,\label{C2E1}\\
&&p(x)=r(x')\ne 0\quad {\rm  and}\quad \lambda \beta(x') = \alpha(x),\label{C2E2}\\
&&p(x)=r(x')=1 \quad {\rm and} \quad \lambda \beta(x') \ne \alpha(x),\label{C2E3} \\
&&p(x)=1\quad {\rm  and}\quad r(x')=0.\label{C2E4}
\end{eqnarray}
\item[(c)] If for every $x\in M$ either \eqref{C2E1} or \eqref{C2E2} holds then $\Phi$ is also a homothety with the dilatation $\Lambda=\lambda$. Moreover, $\varphi$ is 
totally geodesic.
\item[(d)] If \eqref{C2E3} or \eqref{C2E4} holds globally then for every $v,w,Z\in T_x M$
\begin{equation}\label{C2E5} 
\la\Pi(\varphi_\ast v,\varphi_\ast Z), \Pi(\varphi_\ast w,\varphi_\ast Z)\ra'  = C\la v,w\ra \la Z,Z\ra,
\end{equation}
where $C=\lambda (\alpha(x)-\lambda \beta(x'))\ne 0$ in the  case \eqref{C2E3}, and $C=\lambda\alpha(x)\ne 0$ in the case \eqref{C2E4}. In 
particular, $\varphi$ is not totally geodesic.
\end{itemize}
\end{lem}

\begin{proof} 
Assume that $v\ne 0$. Take vectors $A\in T_Z(TM)$ and $A'\in T_{Z'}(TL')$ as in Lemma \ref{CM_nable} and 
Corollary \ref{CM_cor_-1}. Moreover we may assume that $A$ and $A'$ are horizontal with respect to $\nabla$ and $\bar \nabla$, respectively, i.e.,
$K(A)=0$ and $\bar K(A')=0$.

Put $U=U_x$ and $L'=L_{x'}$. Let $J:TL'\to TM'$ be the inclusion map. Put $\bar h=J^\ast h'$.
Since $\varphi: (U, g)\to (L',\bar g)$ and $\Phi:(TU,h)\to (TL,\bar h)$ are conformal diffeomorphisms, so are $\varphi^{-1}:(L',\bar g)\to (U, g)$ and $
\Phi^{-1}:(TL',\bar h)\to (TU, h)$. Dilatations of $\varphi^{-1}$ and $\Phi^{-1}$ are equal to 
$\mu = (1/\lambda)\circ \varphi^{-1}$ and
$\hat \mu = (1/\Lambda)\circ \Phi^{-1}$, respectively.
Thus we have
\begin{eqnarray*}
\|\Phi_\ast A\|'^2 &=& \Lambda(Z)\|A\|^2,\\
\|\Phi^{-1}_\ast A'\|^2 &=& \hat \mu(Z')\|A'\|'^2.
\end{eqnarray*}

Put  for a while $S=S(v,Z)$, $S'=S'(v',Z')$ and $\Pi'=\Pi(v',Z')$.
 Since $\pi_\ast A= v$, $\pi'_\ast A'=v'$, $K(A)=0$, $K'(A')=\bar K(A')+\Pi'=\Pi'$ and 
$Z'$ is orthogonal to $\Pi'$, we have
\begin{eqnarray*}
\|A\|^2 &=& |v|^2,\\
\|A'\|'^2 &=& |v'|'^2 + \omega_\beta(Z')^r|\Pi'|'^2
\end{eqnarray*}
Next applying Lemma \ref{CM_nable} and Corollary \ref{CM_cor_-1}, the equalities
$\pi'_\ast \Phi_\ast A = v'$ and $\pi_\ast \Phi^{-1}_\ast A'=v$, and the fact that $\varphi_\ast S$ is orthogonal to 
$\Pi'$ we get
\begin{eqnarray*}
\|\Phi_\ast A \|'^2 &= &|v'|'^2 + \omega_\beta(Z')^r\big(\lambda(x)|S|^2 + s\lambda^2(x)\la S,Z\ra^2 + |\Pi'|'^2\big)\\ 
\|\Phi^{-1}_\ast A' \|^2 &=& |v|^2 + \omega_\alpha(Z)^p\big(\mu(x') 
|S'|'^2 + q\mu^2(x') \la S',Z'\ra'^2\big)
\end{eqnarray*}

Combining now above equalities and using the definitions of $\mu$ and $\hat \mu$ we get
\begin{eqnarray}
&&\Lambda |v|^2 = |v'|'^2 + \omega_\beta(Z')^r\big(\lambda|S|^2 + \lambda^2\la S,Z\ra^2 + |\Pi'|'^2\big),\label{TL_E1}\\
&&\frac{1}{\Lambda} \big(|v'|'^2 + \omega_\beta(Z')^r|\Pi'|'^2\big) = |v|^2 + \omega_\alpha(Z)^p\big(\frac{1}{\lambda} 
|S'|'^2 + \frac{1}{\lambda^2}  \la S',Z'\ra'^2\big)\label{TL_E2},
\end{eqnarray}
where $\lambda=\lambda(x)$ and $\Lambda=\Lambda(Z)$.
Multiplying equations \eqref{TL_E1} and \eqref{TL_E2}  
side by side we conclude that
\[0 = \omega_\beta(Z')^r\lambda|v|^2|S|^2 + \text{non-negative expression}. \]
Since $v\ne 0$, $S(v,Z)=0$. Since $x\in M$, $v,Z\in T_x M$ were arbitrary the tensor field $S$ vanishes identically. Therefore, $\lambda$ is a constant
function and thus $\varphi$ is a homothety. Hence (a) is proved.

Substituting $S=0$ in \eqref{TL_E1} we get
\[ |\Pi(v',Z')|'^2 = \frac{\Lambda(Z)-\lambda}{ \omega_\beta(Z')^r} |v|^2.\]
Applying Lemma \ref{C1} we get
\begin{equation}\label{C2EP}
|\Pi(\varphi_\ast v,\varphi_\ast Z)|'^2 = \lambda ((1+\alpha|Z|^2)^p-(1+\lambda \beta |Z|^2)^r)|v|^2.
\end{equation}
Using the facts that 
the map $(v,Z)\mapsto |\Pi(\varphi_\ast v,\varphi_\ast Z)|'^2$ is non negative and symmetric with respect to $v,Z$, we conclude (b).

If \eqref{C2E1} or \eqref{C2E2} holds then by \eqref{C11} it follows that $\Lambda=\lambda$. Moreover, in these cases 
\eqref{C2EP} becomes $|\Pi(\varphi_\ast v,\varphi_\ast Z)|'^2=0$. This proves (C).

If \eqref{C2E3} or \eqref{C2E4} holds then it is an elementary computation to check that $\Pi$ satisfies \eqref{C2E5}, proving (d). 
\end{proof}

\begin{lem}\label{Lem_CC2} Suppose that $\dim M\geq 3$ or $\dim M'\leq \dim M+ 1$. Then under the assumptions of Lemma \ref{C2} we have:
$\varphi $ is totally geodesic homothety, $\Phi $ is a homothety and its dilatation $\Lambda$ is equal to $\lambda$.
\end{lem}

\begin{proof}
It suffices to show that under the assumptions the conditions \eqref{C2E3} and \eqref{C2E4} cannot hold. Then the assertion follows from 
Lemma \ref{C2} (a) and (c).
Suppose that \eqref{C2E3} or \eqref{C2E4} holds. Then by Lemma \ref{C2} (d) it follows that the symmetric 
 bilinear form 
$B:T_{x}M \times T_x M \to T_{x'}M'$ given by
$B(v,w)=\Pi(\varphi_\ast v,\varphi_\ast w)$ satisfies the condition \eqref{SS_AL_E1} with $C\ne 0$. If $\dim M\geq 3$ then we have a 
contradiction with Lemma \ref{Lem_AL}, if $\dim M'\leq \dim M +1$ then we have a 
contradiction with Lemma \ref{SS_AL_PRP}.
\end{proof}

\subsection{Main results}  We begin with some definitions. Suppose $\bar M$ is a submanifold of a Riemannian manifold $(M',g')$. 
Suppose that a real-valued non-negative function $C$ on $\bar M$ is given.
We say that 
$\bar M$ is {\em optimal with a coefficient $C$} if for every $x'\in \bar M$ 
the second fundamental form $\Pi$ of $\bar M$ at $x'$  satisfies \eqref{SS_AL_E1} with the constant $C(x')$ that is
\[
\la\Pi(u,w),\Pi(v,w)\ra=C(x')\la u, v\ra\la w,w\ra,\quad u,v,w\in T_{x'}\bar M.
\]

In particular, every totally geodesic submanifold is optimal with the coefficient $0$. By Lemma \ref{Lem_AL} and Lemma \ref{SS_AL_PRP}
it follows that if $\dim \bar M\geq 3$ or ${\rm codim\,} \bar M\leq 1$ then each optimal submanifold is totally geodesic. 

\begin{rmk} \label{MR_R1}
Observe that if $\varphi:M\to \bar M$ is a conformal diffeomorphism such that \eqref{C2E5} holds then 
$\bar M$ is optimal with the coefficient $C/(\lambda\circ\varphi^{-1})^2$.
\end{rmk}

\begin{prp}\label{PRP_OPT} Suppose $\dim \bar M=2$.
Denote by $\kappa'$ and $\bar \kappa$ the sectional curvatures of $M'$ and $\bar M$, and let $\sigma = T_{x'} \bar M$.
If $\bar M$ is optimal submanifold of $M'$ with a coefficient $C$ then $\bar M$ is minimal submanifold of $M'$ and 
$\bar \kappa(\sigma) = \kappa'(\sigma) - 2C(x')$. In particular, if $C$ is constant and $M'$ is a space of constant curvature then so is $\bar M$.  
\end{prp}
\begin{proof}
The fact that $\bar M$ is minimal follows immediately from Lemma \ref{SS_AL_PRP}: it suffices to calculate the trace of the bilinear form 
given by \eqref{AL_E}. The second statement  follows from \eqref{SS_AL_E1}, \eqref{SS_AL_E01} and the Gauss Equation.
\end{proof}

Suppose now that two Riemannian manifolds $(M,g)$ and $(M',g')$ are given and $\dim M\leq \dim M'$. Equip their tangent bundles $TM$ and 
$TM'$ with $(p,q,\alpha)$-metric $h$ and $(r,s,\beta)$-metric $h'$ respectively.
Suppose next that the functions $p,q,r,s,\alpha,\beta$ are constant, and
 $\varphi:M\to M'$ is an imbedding (injective immersion). Let $\bar M=\varphi(M)$. 
 
\begin{thm} \label{MR_T1}
Let $\dim M\geq 3$ or $\dim M'\leq \dim M +1$.
\begin{itemize}
\item[(I)]
Suppose that $\varphi$ is a conformal mapping with a dilatation $\lambda$. Then $\Phi=\varphi_\ast:TM\to TM'$ is conformal if and only if 
$q=\lambda (s\circ\varphi^{-1})$, $\varphi$ is a homothety, $\bar M$ is totally geodesic and for every $x\in M$ $(x'=\varphi(x))$ one of the conditions \eqref{C2E1} or \eqref{C2E2} holds.
\item[(II)]
If $\Phi$ is a conformal mapping then $\varphi$ and $\Phi$ are homotheties and $\Lambda = \lambda$. 
\end{itemize}
\end{thm}

\begin{thm}\label{MR_T2} 
Let $\dim M = 2$ and $\dim M'\geq \dim M +2$.  
\begin{itemize}
\item[(III)]
Suppose that $\varphi$ is a conformal mapping with a dilatation $\lambda$. Then $\Phi=\varphi_\ast:TM\to TM'$ is conformal if and only if $\varphi$ is a homothety,
$q=\lambda (s\circ\varphi^{-1})$, $\bar M$ is optimal with the coefficient
\[
C=\frac{1}{\lambda}((p\alpha)\circ\varphi^{-1}-\lambda \beta r)
\]
and for every $x\in M$ $(x'=\varphi(x))$ one of the properties \eqref{C2E1}-\eqref{C2E4} is satisfied.
\item[(IV)] 
Suppose $\Phi$ is a conformal mapping. Then $\varphi$ a homothety and
\begin{itemize}
\item[(IV1)] If for every $x\in M$ one of the conditions \eqref{C2E1} or \eqref{C2E2} holds then $\Phi$ is also a homothety and $\Lambda = \lambda$. 
\item[(IV2)] If one of the conditions \eqref{C2E3} or \eqref{C2E4} holds globally then $\varphi$ is a minimal
 immersion and for every plane $\sigma=\varphi_\ast (T_x M)$, $x\in M$, the Gauss curvature $\kappa(\sigma)$ of $M$
is 
\begin{equation}\label{MR_T2_cur} 
\kappa(\sigma) =\lambda \kappa'(\sigma) - 2 C(\varphi(x))\lambda, 
\end{equation}
where $\kappa'(\sigma)$ is the Gauss curvature of $M'$. Moreover, $\Phi$ is not a homothety. Its dilatation $\Lambda$ is
\begin{equation}
\Lambda(Z)= \lambda \frac{1+\alpha(x) g(Z,Z)}{1+\lambda \beta(x')r(x') g(Z,Z)}, \quad Z\in T_x M, x'=\varphi(x)\label{MR_T2_E1}.
\end{equation}
\end{itemize}
\end{itemize}
\end{thm}

\begin{proof}[Proof of Theorem \ref{MR_T1} and Theorem \ref{MR_T2}] (I, III $\Rightarrow$) Suppose $\varphi$ and $\Phi$ are conformal mappings. By Lemma \ref{C2}, $\varphi$ is a homothety and for every $x\in M$ ($x'=\varphi(x)$) one of the conditions \eqref{C2E1}-\eqref{C2E4} holds. By Lemma \ref{C1}, $q=\lambda (s\circ\varphi^{-1})$. If $\dim M\geq 3$ or $\dim M'\leq \dim M +1$, by Lemma \ref{C2} (d) and Lemma \ref{Lem_CC2} conditons \eqref{C2E3} and \eqref{C2E4} cannot hold. Therefore Lemma \ref{C2} (c) implies that $\bar M$ is totally geodesic. Moreover, if \eqref{C2E3} or \eqref{C2E4} is satisfied then by Lemma \ref{C2} (d) and Remark \ref{MR_R1}, $\bar M$ is optimal with the coefficient $C=(1/\lambda)((p\alpha)\circ\varphi^{-1}-\lambda \beta r)$. 

(I, III $\Leftarrow$) Taking horizontal (resp. vertical) $A\in T_ZTM$, computing $h'(\Phi_\ast A,\Phi_\ast A)$ and applying relations between $p,q,r,s,\alpha,\beta$ and $\lambda$ one can conclude that $\Phi$ is conformal. 

(II) Suppose $\Phi$ is conformal. By Lemma \ref{C3}, $\varphi$ is also conformal. Therefore (I) and Lemma \ref{C1} imply that $\Phi$ and $\varphi$ are homotheties and $\Lambda=\lambda$. 

(IV) As above we conclude that $\varphi$ is conformal. 

(IV1) It is a consequence of Lemma \ref{C2} (c). 

(IV2) Suppose for every $x\in M$ one of the conditions \eqref{C2E3} or \eqref{C2E4} holds. By Proposition \ref{PRP_OPT}, (I) and the fact that the curvature under the action of a homothety with dilatation $\lambda$ is scaled by $1/\lambda$, $\varphi$ is a minimal immersion and \eqref{MR_T2_cur} holds. Conditions \eqref{C2E3}, \eqref{C2E4} and equation \eqref{C11} imply \eqref{MR_T2_E1}.
\end{proof}

As a direct consequence of of Theorem \ref{MR_T1} we obtain

\begin{cor}\label{Cor_MRT1}
Suppose $\dim M\geq 3$ or $\dim M'\leq \dim M+1$. Let $\varphi:(M,g)\to (M',g')$ be an  imbedding. 
Then we have:
\begin{itemize}
\item[(a)] $\Phi:(TM,h_S)\to (TM',h'_S)$ is conformal if and only if $\varphi$ is totally geodesic homothety.
\item[(b)] $\Phi:(TM,h_{CG})\to (TM',h'_{CG})$ is
conformal if and only if $\varphi$ is totally geodesic isometric imbedding.
\item[(c)] $\Phi:(TM,h_{CG})\to (TM',h'_S)$ is never conformal.
\end{itemize}
\end{cor}

\subsection{An example to Theorem \ref{MR_T2}}\label{SS_EX} 
It is important to show that there is essential difference between Theorems \ref{MR_T1} and \ref{MR_T2}. To do this we give an example of $2$-dimensional manifold $M$, $4$-dimensional manifold $M'$ and an immersion 
$\varphi:M\to M'$ such that $\bar M= \varphi (M)$ is optimal but not totally geodesic. 

Let $\Sigma^d(\rho)$ denote Euclidean $d$-dimensional sphere of radius $\rho$ centred at the origin in $\mathbb{R}^{d+1}$.
Recall
(see \cite{Ch}, Chapter 4 \S 5 page 139) that the {\em second
standard immersion} of $\Sigma^2(1)$ it is a map $\varphi:\Sigma^2(1)\to \Sigma^4(1/\sqrt{3})$ defined as follows:
Consider harmonic homogeneous polynomials $u_i$, $i=1,\dots,5$, in $\mathbb{R}^3$ given by
\begin{eqnarray*}
u_1&=&x_2 x_3, \quad u_2=x_1x_3,\quad u_3=x_1x_2\\
u_4&=& \frac12(x_1^2-x_2^2),\quad u_5=\frac{\sqrt 3}{6}(x_1^2+x_2^2-2x_3^2)
\end{eqnarray*}
and let $u=(u_1,\dots,u_5)$. We define $\varphi$ to be the restriction $u|\Sigma^2(1)$. Then $\varphi:\Sigma^2(1)\to \Sigma^4(1/\sqrt{3})$
is an isometric immersion (but not imbedding). Nevertheless, $\bar M = \varphi (\Sigma^2(1))$ is a minimal submanifold of $\Sigma^4(1/\sqrt{3})$. We show that $\bar M$ is optimal with the constant coefficient $C=1$.

\begin{lem}\label{Lem_EXS}
 Suppose $(M,g)$ is a Riemannian manifold and $u:M\to \mathbb{R}^{d+1}$. Assume that the image $\bar M=u(M)$ is contained in $\Sigma=\Sigma^d(\rho)$
and $u: M\to \Sigma$ is an imbedding. Denote by $\Pi$ and $\bar \Pi$ the second fundamental form of
$\bar M$ in $\Sigma$ and $\bar M$ in $\mathbb{R}^{d+1}$, respectively. Then for every $x'\in \Sigma$ and for every  basis $(e_i)$ of $T_{x'} \bar M$
$$ \la \Pi(e_i,e_k),\Pi(e_j,e_l)\ra =  \la \bar \Pi(e_i,e_k),\bar \Pi(e_j,e_l)\ra - \frac{1}{\rho^2}\la e_i,e_k\ra\la e_j, e_l \ra,$$
where $\la,\ra$ is the canonical inner product in $\mathbb{R}^{d+1}$.
\end{lem}
\begin{proof}
Elementary exercise.
\end{proof}

\begin{prp} If $\varphi$ is the second standard immersion then the submanifold $\bar M=\varphi(\Sigma^2(1))$ is optimal with the constant 
coefficient $C=1$.
\end{prp}
\begin{proof}
Let $M=\Sigma^2(1)$ and $\Sigma= \Sigma^4(1/\sqrt{3})$. Adopt the notations from Lemma \ref{Lem_EXS}. Denote by
$\tilde\varphi$ the restriction of $\varphi$ to lower half sphere $\Sigma^2_-(1)$. 
Since $\varphi(\Sigma^2(1))$ coincides with the closure of $\tilde \varphi(\Sigma^2_-(1))$, it suffices to prove that
$\bar M_-=\tilde \varphi(\Sigma^2_-(1))$ is optimal with the coefficient one.

Let $f$ be the stereographic projection $\Sigma^2(1)\to \mathbb{R}^2$ from the north pole.
Put $\psi = f\circ \tilde\varphi^{-1}$.
Fix $x'=\psi^{-1}(t)$, where $t=(t_1;t_2) \in \mathbb{R}^2$, $|t|<1$. 
Let $e_i=(\partial / \partial 
\psi_i)(x')$.
Put $t^2=t_1^2+t_2^2$ and $t^4 = (t^2)^2$. 
Since $\varphi$ is an isometric immersion and $f$ is the stereographic projection we have 
$$ 
\la e_i,e_j\ra =
 \frac{4\delta_{ij}}{(t^2+1)^2},\quad i,j=1,2,$$ 
where $\delta_{ij}$ is the Kronecker symbol. 
In the light of Lemma \ref{Lem_EXS} (with $\rho = 1/\sqrt 3$), to finish the proof it suffices to show that 
\begin{eqnarray}
\la \bar \Pi(e_i,e_k), \bar \Pi(e_j,e_k) \ra &=& (3\delta_{ik}+1)\frac{16 \delta_{ij}}{(t^2+1)^4}\label{SS_EX_PRPE1}\\
& =&(3\delta_{ik}+1) \la e_i,e_j\ra^2,\quad i,j,k = 1,2.\nonumber 
\end{eqnarray}

 After elementary but laborious calculations we get
\begin{multline*}
\bar \Pi (e_1,e_1)=\frac{4}{(t^2+1)^4}\Big(4t_2(1-t^2-2t_1^2);8t_1(1-t_1^2);4t_1t_2(t_1^2-t_2^2-3);\\
t^4-8t_1^2t_2^2+6t_2^2-6t_1^2+1;\sqrt{3}(t^4-2t^2-4t_1^2+1)\Big)
\end{multline*}
and
\begin{multline*}
\bar \Pi (e_2,e_2)=\frac{4}{(t^2+1)^4} \Big( 8t_2(1-t^2_2);4t_1(1-t^2-2t_2^2);4t_1t_2(t_2^2-t_1^2-3);\\
6t_2^2-6t_1^2 +8t_1^2t_2^2-t^4-1;\sqrt{3}(t^4-2t^2-4t_2^2+1)\Big)
\end{multline*}
and
\begin{multline*}
\bar \Pi (e_1,e_2)=\frac{4}{(t^2+1)^4}\Big( 2t_1(t^2-4t_2^2+1);2t_2( t^2-4t_1^2+1);\\
8t_1^2t_2^2-t^4+1;
4t_1t_2(t_1^2-t_2^2);-4\sqrt{3}t_1t_2\Big).
\end{multline*}
Now one can check that \eqref{SS_EX_PRPE1} holds.
\end{proof}

Now let $\mathbb{R}P^2$ denote the real $2$-dimensional projective space. 
We treat  $\mathbb{R}P^2$ as a Riemannian manifold whose metric $g$ is given by  
the standard two-sheeted covering map  
$\hat \pi:\Sigma^2(1)\to \mathbb{R}P^2$. Put $\hat \varphi(\hat x) = \varphi(x)$ if $\hat x = \hat \pi(x)$. Since $\varphi(x)=\varphi(-x)$, the map $\hat \varphi$ is well defined. Moreover, $\hat \varphi :\mathbb{R}P^2 \to \Sigma^4(1/\sqrt 3)$ is an imbedding.
It is called {\em the first standard imbedding } of  $\mathbb{R}P^2$ into $\Sigma^4(1/\sqrt 3)$. 

Take constants $q,\alpha >0$.
Suppose Cheeger-Gromoll type metrics $h$ and $h'$ on $T(\mathbb{R}P^2)$ and $T(\Sigma^4(1/\sqrt 3))$ are given. 
If

(1) $h=h_{1,q,\alpha+1}$ and 
$h'=h'_{1,q,\alpha}$, or 

(2) $h=h_{1,q,1}$ and $h'=h'_{0,q,1}$\\ 
then
$\hat \varphi_\ast$ is a conformal mapping, but not a homothety. Its dilatation is $\Lambda(Z)= (1+ (\alpha+1)g(Z,Z))/(1 + \alpha g(Z,Z))$ in the case of (1) and 
$\Lambda(Z) = 1+g(Z,Z)$ in the case of (2).

\vskip 1cm
Wojciech Koz\l owski: wojciech@math.uni.lodz.pl\\
Kamil Niedzia\l omski: kamiln@math.uni.lodz.pl\\
Faculty of Mathematics and Computer Science\\
\L \'od\'z University\\
Banacha 22, 90-238 \L \'od\'z\\
Poland

\begin{thebibliography}{99}
\bibitem{BLW} M. Benyounes, E. Loubeau, C. M. Wood, {\em Harmonic sections of Riemannian bundles and metrics of Cheeger-Gromoll type}, Diff. Geom. Appl. {\bf 25} (2007), 322-334.
\bibitem{BLW1} M. Benyounes, E. Loubeau, C. M. Wood, {\em The geometry of generalized Cheeger-Gromoll metrics}, Preprint (arXiv: math.DG/0703059v1).
\bibitem{Ch} B. Chen, {\em Total mean curvature and submanifolds of finite type}, World Scientific, Singapore (1984). 
\bibitem{GK} S. Gudmundsson, E. Kappos, {\em On the geometry of the tangent bundle with the Cheeger-Gromoll metric}, Tokyo J. Math. {\bf 25} (2002), 75-83.
\bibitem{Kw} O. Kowalski, {\em Curvature of the induced Riemannian metric on the tangent bundle of a Riemannian manifold}, J. reine anegew. Math. {\bf 250} (1971) pp. 124-129.
\bibitem{KW} W. Koz\l owski, Sz. M. Walczak, {\em Collapse of unit horizontal bundles equipped with a  metric of Cheeger-Gromoll type}, Diff. Geom. Appl., to appear.
\bibitem{M} M. I. Munteanu, {\em Some aspects on the geometry of the tangent bundles and tangent sphere bundles of a Riemannian manifold}, Mediterr. J. Math.  
{\bf 5} (2008), no. 1, 43--59.
\end{thebibliography}
\end{document}